\documentclass[12pt]{article}
\usepackage[margin=1.6in,top=1.5in,bottom=1.5in]{geometry}
\usepackage{float}
\usepackage{epsfig}
\usepackage{epsf}
\usepackage{here}
\usepackage{color}
\usepackage{pifont}
\usepackage{multirow}

\newcommand{\dpg}[2]{\frac{\partial #1}{\partial #2}}
\newcommand{\ddpg}[2]{\frac{\partial^2 #1}{\partial #2^2}}
\newcommand{\ddp}[3]{\frac{\partial^2 #1}{{\partial #2} {\partial #3}}}

\newcommand{\sine}{\sum _{i=1}^{n}}

\newcommand{\LOi}{{l(Y_i;\theta,b_i)}}
\newcommand{\hlm}{{\rm hl}}
\newcommand{\LOy}{{l(y;\theta,b_i)}}
\newcommand{\LOia}{{l}^{a_i}_{i}}
\newcommand{\hLOia}{{\rm hl}_i}

\newcommand{\Eet}{{\rm E}_{P^*}}
\newcommand{\var}{{\rm var}}

\newcommand{\zi}{\mbox{\boldmath $z_l ^{i}$}}
 
\newcommand{\bi}{b_{i}}

\newcommand{\bet}{\mbox{\boldmath $\beta_l$}}

\newcommand{\X}{\mbox{\boldmath $X$}}

\newcommand{\pl}{\phi_l}
\newcommand{\xiil}{\tilde \xi_{l}^{i}}

\newcommand{\xii}{\mbox{\boldmath $\xi$}^{i}}

\newcommand{\btau}{\mbox{\boldmath $\tau$}}

\newcommand{\tetau}{\hat \theta^{\mbox{\boldmath $\tau$}}}
\newcommand{\phl}{{\rm PHL}_n(\theta)}
\newcommand{\hl}[1]{{\rm HL}(\theta,#1)}
\newcommand{\hls}{{\rm HL}}

\def\inprob{\rightarrow_p}
\newcommand{\bb}{\mbox{\boldmath $b$}}
\newcommand{\bbh}{\mbox{\boldmath $\hat b$}}
\newcommand{\bhl}{\mbox{\boldmath $\hat b$}(\theta)}
\newcommand{\bhltt}{\mbox{\boldmath $\hat b$}(\hat \theta^{\mbox{\boldmath $\tau$}})}
\newcommand{\bihl}{\hat b_i(\theta)}
\newcommand{\bih}{\hat b_i(\theta)}
\newcommand{\bihr}{\hat b^i_r(\theta)}

\newcommand{\ddlx}{\frac{\partial^2 l_i}{\partial x^2}|_{\theta, \hat b_i(\theta)}}
\newcommand{\ddlz}{\frac{\partial^2 l_i}{\partial z^2}|_{\theta, \hat b_i(\theta)}}
\newcommand{\ddJz}{\frac{\partial^2 J}{\partial z^2}|_{\hat b_i(\theta)}}
\newcommand{\ddlzx}{\frac{\partial^2 l_i}{\partial z \partial x}|_{\theta, \hat b_i(\theta)}}
\newcommand{\ddlxz}{\frac{\partial^2 l_i}{\partial x \partial z}|_{\theta, \hat b_i(\theta)}}
\newcommand{\dbhl}{\frac{\partial \hat b_i}{\partial \theta}|_{\theta}}
\newcommand{\dsp}{\displaystyle}
\newcommand{\bq}{\begin{equation}}
\newcommand{\eq}{\end{equation}}

\newtheorem{Theorem}{Theorem}
\newtheorem{Corollary}{Corollary}
\newtheorem{Lemma}{Lemma}

\begin{document}

%Dan1Jos2Com
\title{Inference in HIV dynamics models via hierarchical likelihood}

\author{D. Commenges $^{1,2}$, D. Jolly $^{1,3}$, H. Putter $^{4}$, \\ and R. Thi\'ebaut$^{1,2}$}
 \maketitle

1 INSERM, Epidemiology and Biostatistics Research Center, Bordeaux, F-33076, France; 
2 University of Bordeaux 2, ISPED, Bordeaux, F33076, France; 
3 University of Bordeaux,IMB, Bordeaux, F33405, France; 
4 University of Leiden, Department of Medical Statistics and Bioinformatics, The Netherlands;
E-mail: daniel.commenges@isped.u-bordeaux2.fr

\maketitle

\vspace{\baselineskip} {\bf Abstract:} HIV dynamical models are often based on non-linear systems of ordinary differential equations (ODE), which do not have analytical solution. Introducing random effects in such models leads to very challenging non-linear mixed-effects models. To avoid the numerical computation of multiple integrals involved in the likelihood, we propose a hierarchical likelihood (h-likelihood) approach, treated in the spirit of a penalized likelihood. We give the asymptotic distribution of the maximum h-likelihood estimators (MHLE) for fixed effects, a result that may be relevant in a more general setting. The MHLE are slightly biased but the bias can be made negligible by using a parametric bootstrap procedure. We propose an efficient algorithm for maximizing the h-likelihood. A simulation study, based on a classical HIV dynamical model, confirms the good properties of the MHLE. We apply it to the analysis of a clinical trial.

\vspace{3ex} Keywords:  algorithm, asymptotic, differential equations, h-likelihood, HIV dynamics models, non-linear mixed effects model, penalized likelihood. \vspace{10mm}

\section{INTRODUCTION}
Since the influential paper of Ho et al. (1995) there has been a strong impetus to develop mathematical models for better  understanding the interaction between HIV and the immune system; see Nowak and May (2000). However the statistical inference in these models has raised major challenges coming from the intrication of identifiability and numerical problems. The first problem is numerical: in general the trajectories of the interesting quantities (e.g. viral load or CD4 counts) are solutions of non-linear differential equations that do not have analytical solutions. The second problem is the identifiability problem: the observations recorded on one subject are not informative enough to estimate all the parameters of the model. The first problem is either avoided, simplifying the models to obtain analytical solutions (Wu and Ding, 1999), or solved by using numerical solvers of ordinary differential equations (ODE); Ramsay et al. (2007) proposed an original approach but did not apply it to a random effect model. The second problem is partly treated by considering that the particular values of the parameters for each subject are realizations of random variables with a given distribution in the population. This puts the problem in the framework of non-linear mixed effects models. Laplace approximation of the numerical integrals involved in the computation of the likelihood has been proposed (Beal and Sheiner, 1982; Lindstrom and Bates, 1990); adaptive Gaussian quadrature is another possibility (see Davidian and Giltinan, 1995). We refer to Wu (2005) for a review of statistical issues in HIV models. Recently a stochastic approximation EM (SAEM) algorithm has been proposed (Kuhn and Lavielle, 2005; Donnet and Samson, 2007). In the specific case of HIV dynamics models a Bayesian approach has been proposed by Putter et al. (2002) and Huang, Liu and Wu (2006), while a special algorithm for computing the likelihood and maximizing using a Newton-like method has been proposed by Guedj, Thi\'ebaut and Commenges (2007). However all these methods present difficulties and can be time-consuming.

The hierarchical likelihood (h-likelihood) has been proposed for generalized linear models with random effects by Lee and Nelder (1996) and further studied in Lee and Nelder (2001) and Lee, Pawitan and Nelder (2006) and for non-linear mixed effects models by Noh and Lee (2008). This is very similar to an approach called penalized likelihood used by McGilChrist and Aisbett (1991) and Therneau and Grambsch (2000) for frailty models. The main idea is to treat the random effects (or the frailties) as parameters and to find estimates of all the parameters by maximizing a function which is essentially the loglikelihood conditional on the random effects minus a penalty term which takes large values if the ``random'' parameters are very dispersed. Penalized likelihood has also been used for function estimation (O'Sullivan; 1988). The advantage of this approach is that it may avoid computing numerical integrals. The curse of dimensionality is transferred from the dimension of numerical integrals to the dimension of the space on which maximization takes place. There are problems with this approach. One is the asymptotic distribution of the estimators of the fixed parameters; another is the estimation of the variances of the random parameters. Consistency of the maximum h-likelihood estimators (MHLE) has not been proved. It is often suggested to revert to the likelihood to have consistent estimators of the fixed parameters, but then the most important benefits of h-likelihood in terms of computational burden is lost. Last but not least is the problem of maximizing a complicated function over several hundred parameters.

The aim of this paper is to develop a (partly non-standard) h-likelihood approach to HIV dynamics models which completely avoids computation of the likelihood. This is in the spirit of penalized likelihood in the sense that we do not try to precisely estimate the variances of the random effects. One aim is to study the asymptotic distribution of the MHLE for a given choice of the penalty. Another aim is to find an efficient maximization algorithm.

The paper is organized as follows. In section 2 we describe a statistical model based on an ODE system in a general form and in a particular form which will be used for simulations. In section 3 we describe h-likelihood and we give the asymptotic distribution of the MHLE for fixed effects when the number of subjects tends toward infinity. We propose a parametric bootstrap procedure to correct the bias of the MHLE. In section 4 we propose a strategy for choosing the penalty based on the guess of an upper bound of the variance of the random effects. An efficient maximization algorithm is presented in section 5. Section 6 presents a simulation study. Section 7 presents the analysis of a clinical trial. We conclude in section 8.

\section{A POPULATION DYNAMICS MODEL}\label{DynMod}
\subsection{A general model for the system}\label{syst}
The dynamics of the concentrations of virions and CD4+ T-cells (in short, CD4) in different stages (represented by ${\X}^{i}(t)$)  can be described by an ODE system. We allow the values of the parameters to vary between subjects; thus we 
consider a population model, as in Guedj, Thi\'ebaut and Commenges (2007). For
subject $i$ with $i=1,...n$, this can be written:
\begin{equation}\label{sys2}
\left \lbrace \begin{array}{l}
\frac{d{\X^{i}}(t)}{dt} = f({\X}^{i}(t), \xii) \\
{\X} ^{i}\small{} (0)=h(\xii)
 \end{array} \right.
\end{equation}
where $\X^{i} (t)= (X^{i} _{1}(t),...,X^{i} _{K}(t))'$ is the
vector of the $K$ state variables (or components);
$\xii=(\xi_{1}^{i},...,\xi_{p}^{i})$ is a vector
of  $p$ individual parameters which appear naturally in the ODE
system and have generally a biological interpretation.
Similarly to generalized (mixed) linear models, we introduce a
link function which relates $\xii$ to a linear
model involving explanatory variables and random effects:
\begin{equation}\label{statist}
\Psi _l (\xi_{l}^{i})=\xiil=\left \{ \begin{array}{ll} \pl + b^i_l+\zi(t) \bet ,& ~~ l=1,\ldots, R,\\
    \pl +\zi(t) \bet ,&~~ l=R+1,\ldots, p,
\end{array}
\right.
\end{equation}
where $\pl$ is the intercept, $\zi(t)$ are vectors of explanatory variables associated with the fixed effects of the $l$th biological parameter; these explanatory variables may be time-dependent, in which case the ODE system has time-dependent parameters. The $\bet$'s are vectors of regression coefficients;  $\bi=(b^i_1,\ldots,b^i_R)$ is the individual vector of random effects. We assume $\bi \sim \mathcal{N}(0,{\mbox{\boldmath $\Sigma$}})$ with $\Sigma$ diagonal with diagonal elements $\tau_l^2$. More general models could of course be considered.

\subsection{Model for the observations}
 Let $Y_{ijm}$ denote the $j$th measurement of the
$m$th observable component for subject $i$ at time $t_{ijm}$; we assume that:
\begin{equation}\label{observation}Y_{ijm}= g_{m}(\X^i(t_{ijm})) + \epsilon_{ijm}, ~~~ i=1,...,n,~~~ j=1,...,n_{im},~~~\end{equation}
for $~ m=1,...,M$, where $g_{m}(.)$ are known functions and
where the $\epsilon_{ijm}$ are independent Gaussian variables with zero mean and variances $\sigma_{m}^{2}$. {The $\epsilon_{ijm}$'s are supposed independent because they represent measurement errors.
The model for the observations may be complicated by the detection limits of assays
leading to left-censored observations $Y_{ijm}$.

\subsection{A particular model for HIV dynamics}\label{standard}
For illustrating the proposed method we present a version of a rather standard model for the HIV dynamics model, close to that used by Nowak and Bangham (1996):

\begin{eqnarray*}\label{ode}
%\frac{dQ^i}{dt}&=& \lambda^i + \rho^i T^i - \alpha^i Q^i- \mu^i_{Q} Q^i \\
\frac{dT^i}{dt}&=& \lambda^i   - \gamma^i T^iV^i - \mu^i_{T}T^i  \\
\frac{dT^{*i}}{dt}&=& \gamma^i  T^iV^i - \mu^i _{T^{*}} T^{*i} \\
\frac{dV^i}{dt}& =& \pi^i T^{*i} - \mu^i _{V}V^i \\
\end{eqnarray*}
where $T^i$, $T^{*i}$ represent the concentrations (implicitly depending on $t$) of
non-infected and  infected CD4 respectively,  and $V^i$ stands for the concentration of virus.

Here the components of $\xii_l=(\lambda^i, \gamma^i, \mu_T^i, \mu_{T^*}^i, \pi^i,\mu_V^i)$ represent rates of events such as production of new cells or particles, rates of infection after meeting between different particles.  As for the $\Psi_l(.)$, we will take the natural log-transform for all the parameters: the natural log-transform can be justified if we think of the parameters as expectations of Poisson variables and has the advantage that the standard deviations of the transformed parameters may be interpreted as coefficients of variations of the estimators of the natural parameters.
 For the simulations, in the link equation (\ref{statist}) we will take random effects  for $\lambda^i$, $\pi^i$ and $\mu^i _{T^{*}}$ (so $R=3$) with the $b^i_l$ being normal and independent with variances $\tau^2_{\lambda}$, $\tau^2_{\pi}$, $\tau^2_{\mu_{T^*}}$, respectively. For all parameters except for $\gamma^i$ we will take no explanatory variable. The effect of the treatment will be modeled as modifying $\gamma^i$ according to the equation:
$$\tilde  \gamma^i=\log \gamma^i=\gamma_0+\beta_1z^i_1(t)+\beta_2z^i_2(t),$$
where $z^i_1(t)$ and $z^i_2(t)$ are treatment indicators. The treatment may change with time; here we will suppose that they are fixed for $t\ge 0$ but take the value $0$ for $t<0$. We assume that at $t=0$ the patients are at the equilibrium of the system with $z^i_1(t)=z^i_2(t)=0$ and this gives important information.
As for the observation equation (\ref{observation}) we will take in the simulations:

\begin{eqnarray*}\label{Obs}
Y_{ij1}&= & \log_{10} V^{i}(t_{ij1}) + \epsilon_{ij1}\\
Y_{ij2}&= & [T^i(t_{ij2})+T^{*i}(t_{ij2})]^{1/4} + \epsilon_{ij2}  \\
Y_{ij3}&= & [T^{*i}(t_{ij3})]^{1/4} + \epsilon_{ij3}
\end{eqnarray*}

\section{THE HIERARCHICAL OR PENALIZED LIKELIHOOD}
\subsection{Asymptotic Distribution of the MHLE}

Let us consider the following model: conditionally on $b_i$, $Y_i$ has a density $f_Y(.;\theta, b_i)$, where $\theta$ is a vector of fixed parameters of dimension $q$ ($\theta \in \Theta \subset \Re ^q$) and $b_i$ are random effects (or parameters) of dimension $R$. The $(Y_i,b_i)$ are independently identically distributed (iid). Typically $Y_i$ is multivariate of dimension $n_i$. We assume that the $b_i$ have density $f_b(.; \btau)$ with zero expectation and where $\btau$ is a vector of parameters. We denote by $P^*$ the true probability and $\theta^*$ and $\btau^*$ the parameter values which specify the distribution of the observed $Y_i$. Typically $Y_i$ is (at least partially) observed while $b_i$ is not.

 Estimators of both $\theta$ and $\bb=(b_1,\ldots,b_n)$ are defined as maximizing the (normalized) extended loglikelihood, called here (by abuse of language) h-loglikelihood:
$$\hl{\bb,\btau}=  L^{\theta, \bb}_{_n}-\frac{1}{n}\sine J(b_i; \btau),$$
where $L^{\theta, \bb}_{n}$ is the loglikelihood (normalized by $\frac{1}{n}$) for the observation conditional on $\bb$, and $J(b_i; \btau)=-\log f_b(b_i; \btau)$.
We denote by $(\hat \theta^{\btau},\bbh ^{\btau})$ the values which maximize $\hl{\bb,\btau}$ for given $\btau$; $\hat \theta^{\btau}$ will be called the MHLE of the parameters $\theta$.
We have $\hl{\bb,\btau}=\frac{1}{n}\sum _{i=1}^n \hlm(Y_i;\theta,b_i, \btau)$ with $\hlm(Y_i;\theta,b_i,\btau)=\LOi-J(b_i; \btau)$, where $\LOi$ is the loglikelihood for subject $i$ conditional on $b_i$. 
For simpler notation we will not always make the dependence in $\btau$ explicit and will write for instance $\hl{\bb}$ for $\hl{\bb,\btau}$. We shall make the additional assumptions:

{\bf A1} $\LOy$ and $J(b_i; \btau)$ are continuous and twice-continuously differentiable functions of $\theta$ and $b_i$ for all $y$ and $\btau$;

{\bf A2} $\Eet \LOi$ exists for all $\theta \in \Theta$.

We shall derive asymptotic results for the MHLE of the fixed parameters $\theta$, which do not require that
 $\btau=\btau^*$.

\begin{Lemma} Under assumptions A1 and A2 the MHLE for fixed effects are M-estimators.
\end{Lemma}
{\bf Proof}.
Consider the profile h-loglikelihood $\phl=\hl{\bhl}$, where $\bhl={\rm argmax}_b~\hl{\bb}$.
$(\tetau, \bhltt)$ maximizes $\phl$, thus $\tetau$ is the profile h-likelihood estimator.  Remembering that $\hl{\bb}= \frac{1}{n}\sum _{i=1}^n \hlm(Y_i;\theta,b_i)$, it is clear that the components of $\bhl$ are the $\bihl={\rm argmax}_{b_i} [\hlm(Y_i;\theta,b_i)]$. Thus $\phl=\hl{\bhl}=\frac{1}{n}\sum _{i=1}^n \hlm(Y_i;\theta,\bihl)$. It follows that $\tetau$ is a M-estimator because it is clear that $\tetau$ is the maximum of $M_n(\theta)=n^{-1}\sum _{i=1}^n m_{\theta}(Y_i)$, where $m_{\theta}(y)$ is a known measurable function: here  $m_{\theta}(y)=\hlm(y;\theta,b(y;\theta))$ where $b(y;\theta)={\rm argmax}_{b} [\hlm(y;\theta,b)]$ (Van der Vaart, 1998, p 41).

For the convergence result we need the additional assumption:

\noindent{\bf A3} For every sufficiently small ball $U\in \Theta$, $\Eet \sup_{\theta\in U}\hlm(y;\theta, b(y;\theta))< \infty$.

In the convergence theorems of the MHLE we will emphasize the fact that it depends on $n$ by writing $\hat \theta^{\btau}=\hat \theta^{\btau}_n$.
\begin{Theorem}  If $\Theta$ is compact and assumption A1-A3 holds, the MHLE of fixed effects $\hat \theta^{\btau}_n$ converges in probability toward $\theta_0^{\btau}=argmax_{\theta}~\Eet [\hlm(Y_i;\theta,\bihl)]$, for any $\btau$.
\end{Theorem}
{\bf Proof}.
By the law of large numbers $M_n(\theta) \inprob M(\theta)$ where $M(\theta)= \Eet [\hlm(Y_i;\theta,\bihl)]$. Let us call $\theta_0^{\btau}$ the value, that we assume unique, at which $M(\theta)$ attains its maximum. The conditions stated in the Theorem, together with the continuity assumption A1, allow us to apply Wald's consistency proof (van der Vaart, 1998, Theorem 5.14, p48).

\begin{Corollary} The MHLE of the fixed parameter of the statistical model described in section \ref{DynMod} converges in probability toward $\theta_0^{\btau}=argmax_{\theta}~\Eet [\hlm(Y_i;\theta,\bihl)]$.
\end{Corollary}
{\bf Proof}.  In the case of the statistical model of section \ref{DynMod} we have $\hlm(Y_i;\theta,\bihl)=\sum _{m=1}^M[ -\frac{n_i}{2} \log \sigma^2_m - \sum _{j=1}^{n_i}\frac{(Y_{ijm}-\phi_m(t_{ijm};\theta,\bihl))^2}{2\sigma^2_m}]-\sum_{r=1}^R \frac{\bihr^2}{2 \tau^2}$, where $\phi_m(t_{ijm};\theta,b_i)=g_{m}(\X^i(t_{ijm}))$ (where $\X^i(t_{ijm})$ is the solution of the ODE system with parameters $\theta,b_i$). In case where $\sigma^2_m$ are fixed, assumption A3 is trivially satisfied because we can remove the terms involving $\sigma^2_m$ and obtain a function which is bounded by zero. If we include the $\sigma^2_m$ in the parameters that we wish to estimate, assumption A3 is  satisfied since  $\hlm(Y_i;\theta,\bihl)\le \sum _{m=1}^M -\frac{n_i}{2} \log \sigma^2_m - \sum _{j=1}^{n_i}\frac{(Y_{ijm}-\phi_m(t_{ijm};\theta,\bihl))^2}{2\sigma^2_m}\le \sum _{m=1}^M -\frac{n_i}{2} \log \tilde \sigma^2_m - n_i/2$, with $\tilde \sigma^2_m=\frac{\sum _{j=1}^{n_i}(Y_{ijm}-\phi_m(t_{ijm};\theta,\bihl))^2}{n_i}$. It seems reasonable to conjecture that $\Eet [-\log {(Y_{ijm}-\phi_m(t_{ijm};\theta,\bihl))^2}] < \infty$. We can compactify the space by taking $\Theta=\bar \Re^d$. If some parameters take an infinite value, $\hlm(Y_i;\theta,\bihl)$ take either the value $-\infty$ or a finite value.

 Now the problem is to investigate whether $\theta_0^{\btau}$ is equal or close to $\theta^*$. $M(\theta)$ can be considered as an approximation of minus the Kullback-Leibler divergence. This is obtained by replacing the expectation in $\bb$ by the mode. The approximation is exact if $\phi_m$ are linear functions in $\bb$ but this is not true in general. However there is a possibility of reducing the bias (see section 3.3).

The asymptotic normal distribution holds for M-estimators under some regularity conditions. We make use of Theorem 5.23 of van der Vaart (1998) which only requires a Lipshitz condition on $m_{\theta}(y)=\hlm(y;\theta,b(y;\theta))$ that we can establish if the following assumption bearing on $u^{\theta}(y)=\dpg{\hlm(y;\theta,\hat b(y;\theta))}{\theta}$ holds. We shall use $u_i^{\theta}=u^{\theta}(Y_i)$ and will give an alternative expression in formula (\ref{utet}).

\noindent{\bf A4} There is a neighborhood $\Theta_0\subset \Theta$ of  $\theta_0^{\btau} $ such that the function $\dot m(y)=\sup_{\theta\in \Theta_0}u^{\theta}(y)$ has the property: $\Eet\|\dot m(Y_i)\|^2 <\infty$.

\begin{Theorem} \label{asympt} Assume assumptions A1-A4 hold.
Then  $\sqrt n (\hat \theta_n^{\btau}-\theta_0^{\btau})$ is asymptotically normal with zero expectation and variance equal to
$$\Sigma(\theta_0^{\btau})= \{\Eet[H_i^{\theta_0^{\btau}}]\}^{-1} \{\Eet[u_i^{\theta_0^{\btau}}u_i^{\theta_0^{\btau}T}]\} \{\Eet[H_i^{\theta_0^{\btau}}]\}^{-1},$$
where  $H_i^{\theta}=\dpg{u_i^{\theta}}{\theta}$.
\end{Theorem}
{\bf Proof}.
The theorem follows by applying Theorem 5.23 of van der Vaart (1998). In this theorem, the main condition is that there exists a measurable function $\dot m$ with $\Eet \dot  m^2 < \infty$ such that for every $\theta_1$ and  $\theta_2$ in a neighborhood $\Theta_0$ of  $\theta_0 $ we have:
\begin {equation} \label{5.23} |m_{\theta_1}(y)-m_{\theta_2}(y)| \le \dot m(y) \|\theta_1-\theta_2\|.\end{equation}

A Taylor series expansion gives: $m_{\theta_1}-m_{\theta_2}= (\theta_1-\theta_2)^T\dpg{m_{\theta}}{\theta}(\tilde \theta)$, where $\tilde \theta \in \Theta_0$. This yields:
$$|m_{\theta_1}-m_{\theta_2}| \le \|\dpg{m_{\theta}}{\theta}(\tilde \theta)\|  \|\theta_1-\theta_2\|\le \sup_{\theta\in \Theta_0} \|\dpg{m_{\theta}}{\theta}(\theta)\|  \|\theta_1-\theta_2\|.$$
Then assumption A4 allows us applying the Theorem 5.23 of van der Vaart (1998).

 For applying Theorem \ref{asympt}, it remains to compute the first and second derivatives of $m_{\theta}(Y_i)$ in terms of derivatives of the likelihood conditional on the random effects. We write $l_i(\theta,\bih)=l(Y_i;\theta,\bihl)$. Let us call $\dpg{l}{x}$ (resp.$\dpg{l}{z}$) the derivatives of $l_i(.,.)$ wrt the first (resp. the second) argument and $\dpg{J}{z}$ the derivative of $J(.)$ wrt its argument. We have
$u_i^{\theta}=\dpg{l_i}{x}|_{\theta, \bih}+\dpg{l_i}{z}|_{\theta, \bih}\dpg{\hat b_i}{\theta}|_{\theta}-\dpg{J}{z}|_{\bih}.$ However, because $\bih$ maximizes $\hlm(Y_i; \theta,b)$ we have 
\begin{equation} \label{eqbhl}\dpg{l_i}{z}|_{\theta, \bih}\dpg{\hat b_i}{\theta}|_{\theta}-\dpg{J}{z}|_{\bih}=0.\end{equation}
 Hence we obtain that
\begin{equation} \label{utet} u_i^{\theta}=\dpg{l_i}{x}|_{\theta, \bih}.\end{equation}
 That is $u_i^{\theta}$ is simply the derivative of the loglikelihood as if $b$ was fixed, computed in $(\theta,\bih)$.

Next we have $H_i^{\theta}=\dpg{u_i^{\theta}}{\theta}=\ddlx +\ddlzx \dbhl$. Differentiating equation (\ref{eqbhl}) wrt $\theta$ we have:
$$\ddlxz+\ddlz \dbhl -\ddJz \dbhl=0,$$
from which we obtain:
$$\dbhl= -\left [\ddlz-\ddJz\right ]^{-1}\ddlxz.$$
Hence:
$$H_i^{\theta}=\ddlx -\ddlzx \left [\ddlz-\ddJz\right ]^{-1}\ddlxz.$$

 In practice we can plug in the estimator $\hat \theta^{\btau}$ to obtain an estimator of $\Sigma(\theta_0^{\btau})$ (using the continuous mapping theorem). We may also use the observed scores and Hessian. By virtue of the law of large numbers they converge toward their expectations, and again the continuous mapping theorem allows to prove consistency of the resulting estimator.

\subsection{Correction of the bias}\label{bootstrap}
We have shown in section 3.2 that  the MHLE $\hat \theta^{\btau}$ tends toward $\theta_0^{\btau}$ which is in general different from $\theta^*$; thus there is an asymptotic bias $\theta_0^{\btau}-\theta^*$. Note that the asymptotic distribution is valid for any $\btau$, and on the other hand, $\hat \theta^{\btau}$ is biased even for $\btau=\btau^*$. Thus the problem of this approach is essentially that of the bias, although a small bias may be acceptable if it goes with a small variance. We propose to partially correct the bias by parametric bootstrap (Efron and Tibshirani, 1993).

Specifically, for $s=1,\ldots,S$, generate the $b_i^{s}$ from $f_b(.,\btau)$; generate $Y_i^s$ from $f_Y(.;\hat \theta^{\btau}, b_i^{s})$; compute the MHLE $\hat \theta^{\btau, s}$ for these data. An estimator of the bias is $S^{-1} \sum_{s=1}^S(\hat \theta^{\btau, s}-\hat \theta^{\btau})$. Thus the corrected estimator, called cMHLE, is
$$\check \theta^{\btau }=\hat \theta^{\btau }-S^{-1} \sum_{s=1}^S(\hat \theta^{\btau, s}-\hat \theta^{\btau}).$$

This correction slightly increases the variance. The variance of $\check \theta^{\btau }$ can be computed through the formula $\var~ \Eet  (\check \theta^{\btau }|\hat \theta ^{\btau}) + \Eet \var (\check \theta^{\btau }|\hat \theta ^{\btau})$. Neglecting the bias of the MHLE in this computation we obtain:
$$\var~ \check \theta^{\btau}\approx (1+S^{-1})\var~ \hat \theta^{\btau}$$

\section{PENALTY CHOICE}\label{crossvalidation}

Profile likelihood has been proposed by Therneau and Grambsch (2000) and Lee and Nelder (2001) but it has the drawback of requiring the computation of the marginal likelihood. We propose a strategy for penalty choice which avoids this computation. For any choice of $\btau=(\tau_1,\ldots,\tau_R)$ we have that $\hat \theta^{\btau}$ has an asymptotic normal distribution with expectation $\theta_0^{\btau}$ and with a variance that can be estimated. We propose to take a reasonable upper bound of $\btau$, that is, the value $\btau^u=(\tau^u,\ldots,\tau^u)$ where $\tau^u$ is considered as an approximate upper bound for the $\tau^*_i$. First, note that since we are working with natural logarithms of the biological parameters, the $\tau_i$ may be interpreted as coefficients of variation of these parameters. It seems reasonable (and is in agreement with the literature) to expect coefficients of variations of parameters such as rate of production of new lymphocytes ($\lambda$) or death rate of uninfected lymphocytes ($\mu_T$) are not very large, that is no more than $0.3$.

\section{MAXIMIZATION ALGORITHM}

 Newton-like algorithms use an approximation of the Hessian of the function to maximize. Since there are many parameters, this matrix can be very large. For instance in our application $q=7$, $R=3$, $n=100$, so the number of parameters is $q+nR= 307$. In complex problems, both gradient and Hessian have to be computed numerically. Particular care must be spent to compute the Hessian both economically and precisely. The algorithm we propose is an adaptation of the Marquardt algorithm (Marquardt, 1963), taking advantage of the special structure of the Hessian in our problem. We draw two consequences of this special structure: (i) there are many terms which are equal to zero, so we do not need to compute them; (ii) the matrix is not far from being block-diagonal.

We shall  first consider the particular case where the number of random and fixed
effects are equal ($R=q$) and the loglikelihood of subject $i$, $\LOi$, depends only on  $\theta + b_i$.
We are interested in maximizing the following function:
$$ \hl{\bb}=\frac{1}{n}\sum_{i=1}^n \left [\LOi-\sum_{r=1}^R\frac{{b_r^i}^2}{2\tau^2}\right ]\cdot $$
It is useful to reparametrize in term of $a_i = \theta + b_i$. One finds
$$ \hls=\frac{1}{n}\sum_{i=1}^n \left [\LOia-\sum_{r=1}^R\frac{{(a_r^i-\theta_r)}^2}{2\tau^2} \right ]=\frac{1}{n}\sum_{i=1}^n \hlm_i \cdot$$
With this parameterization the loglikelihood, $\LOia=l(Y_i;\theta,a_i-\theta)$, which is the complex part, depends only on $a_i$ so that many derivatives of the h-loglikelihood are very simple:%
\begin{equation} \label{eq1} \frac{\partial \hls}{\partial \theta_r} = \frac{1}{n} \sum_{i=1}^n \frac{a^i_r - \theta_r}{\tau^2};\end{equation}

$$ \frac{\partial^2 \hls}{\partial \theta_r \partial a_{r'}^i} =  \frac{\delta_{rr'}}{n \tau^2}~;~
\frac{\partial^2 \hls}{\partial a_r^i \partial a_{r'}^{i'}} = 0, \textrm{if $i \neq i'$ ; } \frac{\partial^2 \hls}{\partial \theta_r \partial \theta_{r'}} = - \frac{\delta_{rr'}}{\tau^2},
$$
where $\delta_{rr'}=1$ if and only if $r=r'$.
This leads to a specific block structure of the Hessian matrix, involving blocks $A =\ddpg{\hls}{\theta}= -\frac{1}{\tau^2} I_R$ and $D=\ddp{\hls}{\theta}{a_i}= \frac{1}{n \tau^2} I_R$ (where $I_R$ is the identity matrix of dimension $R$; $D$ does not depend on $i$) and $C_i= \frac{1}{n}\ddpg{\hLOia}{a_i}$; the structure is displayed in Figure \ref{hessian1}.

\begin{figure}[H]
\begin{center}
\includegraphics[angle=270,scale=0.4]{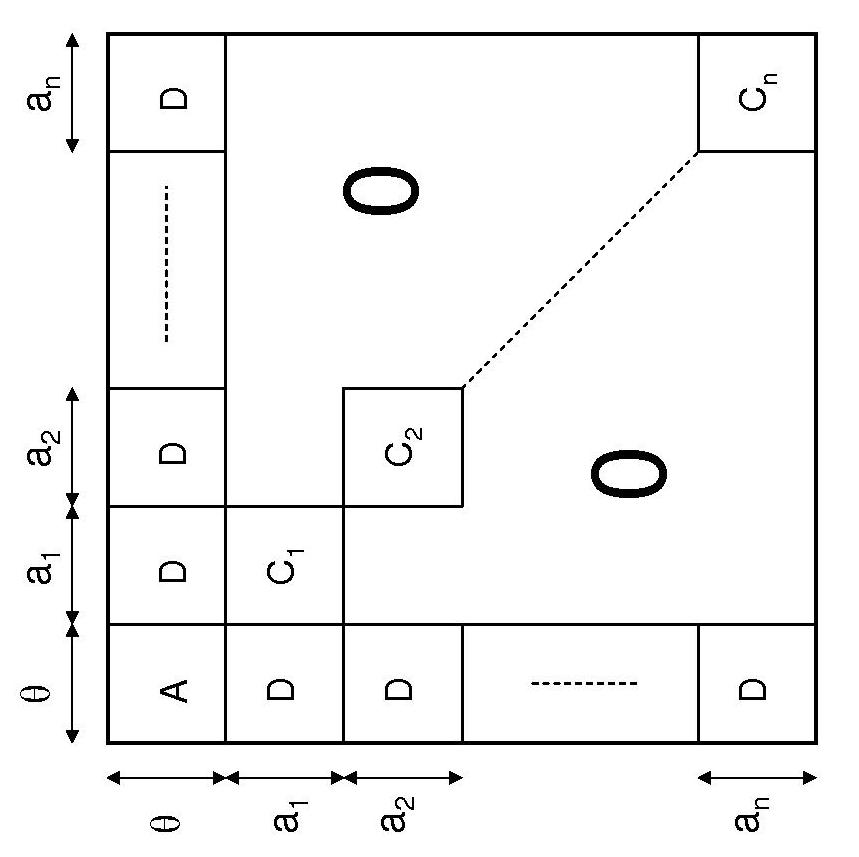}
\caption{\label{hessian1}Hessian matrix in the case $R=q$.
$\dsp A = -\frac{1}{\tau^2} I_R$ and $\dsp D = \frac{1}{n \tau^2}
I_R$, $I_R$ is the identity matrix of dimension $R$, and $C_i= \frac{1}{n}\ddpg{\hLOia}{a_i}$.}
\end{center}
\end{figure}
Fast computation of this large $(n+1)R\times (n+1)R$ Hessian matrix is possible for two reasons: (i) only the terms of blocks $C_i$  require computation of the likelihood; (ii) for computing the terms of block $C_i$ we only need to compute the second derivatives of $\LOia$ (and not of the whole h-likelihood). Finally there are $nR(R+1)/2$ terms to compute, each involving only one computation of the solution of the ODE system (needed for the numerical differentiation): it follows that the number of computations of the solution of ODE system does not exceed that required for the computation of the Hessian for an ordinary (without random effect) non-linear model with $q$ parameters ! 

The nearly diagonal structure of the Hessian led us to design the so-called  ``patient-by-patient'' algorithm, decoupling the optimization between patients. Denoting by $a_i(k)$ and $\theta(k)$ the values at iteration $k$, iteration $k+1$ proceeds in two steps:

Step 1: For $i=1,\ldots,n$:
make one Marquardt step for optimizing the function $\LOia-\sum_{r=1}^R\frac{{(a_r^i-\theta_r(k))}^2}{2\tau^2}$ on $a_i$; this gives $a_i(k+1)$; Step 2:
compute  $\theta_r(k+1) = \frac{1}{n} \sum_{i=1}^n {a}^i_r(k+1) $  (which satisfies (\ref{eq1}));
go to step 1 (until convergence is reached).

The patient-by-patient algorithm works very well far from the maximum when the global Marquardt algorithm is hampered by the need of a large increase of the diagonal of the Hessian. However the decoupling between patients also leads to a loss of efficiency so that close to the maximum it is less efficient than the global Marquardt algorithm. This observation led us to devise a hybrid algorithm:
use the patient-by-patient algorithm until all blocks $C_i$ are definite-positive; then switch to the global Marquardt algorithm. Note that ensuring that all blocks $C_i$ be definite-positive does not imply that the Hessian is so; generally however it is not far from being the case so that the Marquardt algorithm is efficient.

We now consider the case where there are $R$ fixed parameters, that we call $\alpha$, associated with a random effect; such as above the loglikelihood of subject $i$, $\LOi$, depends only on  $\alpha + b_i$ and a vector of fixed parameters $\beta$. As in the preceding case, the Hessian has a particular structure (see Figure \ref{hessian2}). It involves the blocks $A$, $D$ and $C_i$ as above, and in addition blocks $B=\ddpg{\hls}{\beta}$ and $B_i=\frac{1}{n}\ddp{\hlm_i}{\beta}{a_i}$.
\begin{figure}[H]
\begin{center}
\includegraphics[angle=270,scale=0.40]{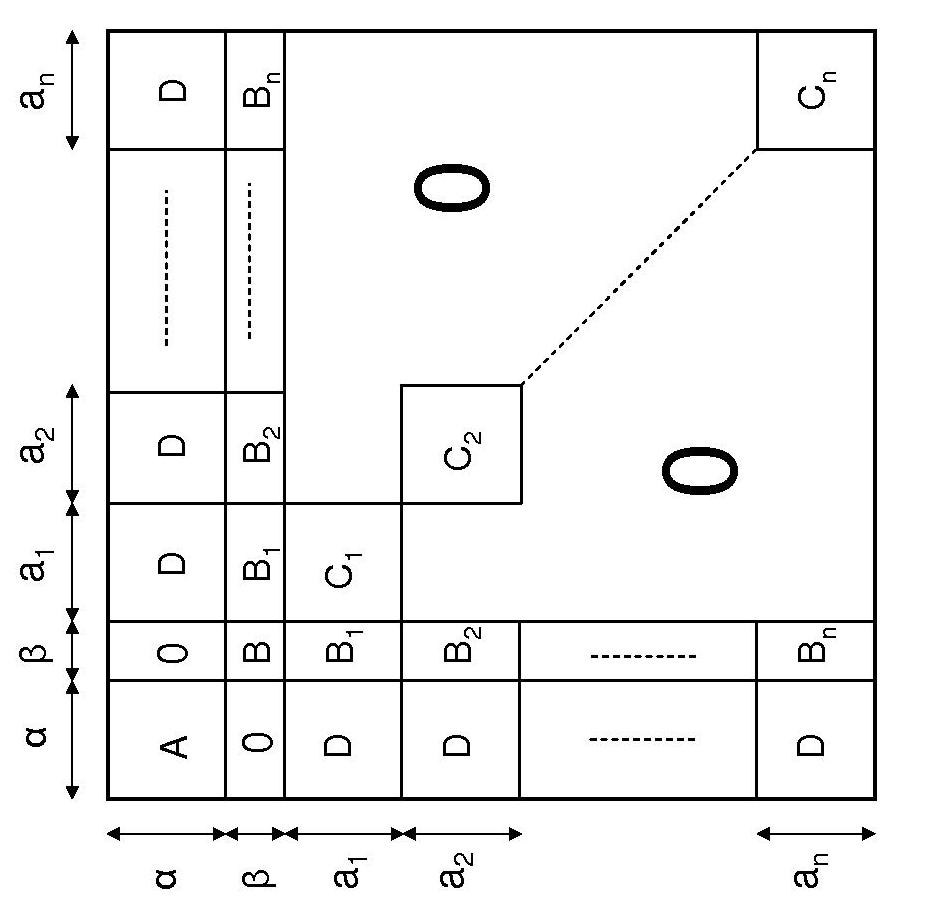}
\end{center}
\caption{\label{hessian2}Hessian matrix in fixed and random effects
case where $\dsp A = -\frac{1}{\tau^2} I_R$, $\dsp D  = \frac{1}{n
\tau^2} I_R$ (where $I_R$ is the identity matrix of dimension $R$), $B=\ddpg{\hls}{\beta}$, $B_i= \frac{1}{n}\ddp{\hlm_i}{\beta}{a_i}$ and $C_i=\frac{1}{n} \ddpg{\hlm_i}{a_i}$}
\end{figure}%

The idea, like previously, is to deal with the case where there are
non-definite positive $C_i$. For that, we use the two steps of the patient-by-patient
approach which give the individual parameters ${a}_i(k+1)$ and
their means ${\alpha}(k+1)$. Then keeping these values fixed, we find the other fixed
parameters ${\beta}(k+1)$ by a step of Marquardt
algorithm with the block $B$. As soon as all blocks $C_i$ and the
block $B$ are definite positive, we switch to the global Marquardt algorithm.

\section{A SIMULATION STUDY}\label{simulation}

\subsection{Description of the simulation study}
We did simulations from the model described in section \ref{standard}. We fixed (that is we did not estimate) the parameters $\tilde \mu_V$, $\tilde \mu_T$ and $\sigma_i$, at values which are plausible in view of the literature (taking as time unit the day and as volume unit the micro-liter):  $\tilde \mu_V = 3.40;~ \tilde \mu_T = -2.20;~ \sigma_i=0.5,~ i=1,2,3$. The values for the other parameters (to be estimated), including the two treatment effects $\beta_1$ and $\beta_2$, are given in Table 2. For each replica, observations for $n=100$ subjects were generated; for each subject  $n_{im}=10$ observations for the three compartments ($m=1,2,3$) were generated at times $0, 3, 6, 9, 12, 15,18,21,24 ,30$.

\subsection{Efficiency of the algorithm}
We did a simulation to compare the number of iterations of the global Marquardt algorithm and the hybrid algorithm. We tried the two algorithms with models including one to three random effects. The initial values were: $\tilde \lambda=5.0;~\tilde \mu_{T^*}=0;~\tilde \pi=0;~\tilde \gamma_0=-5.0;~\beta_1=-1.0;~\beta_2=-1.0$. The global Marquardt algorithm did not always converge in less than 150 iterations while the hybrid algorithm nearly always converged (see Table \ref{meaniter}); when they both converged, this was toward the same values (close to the true parameter values). We checked that when we started from different values the algorithms converged toward the same values. The hybrid algorithm is faster than the global one. In Table \ref{meaniter} we give the mean number of iterations until convergence (computed on 100 replications) for which the algorithm converged in less than 150 iterations. For instance the mean number of iterations for 3 random effects was 25 versus 71 for the hybrid versus the global algorithm. The mean time of one iteration is about the same for the two algorithms. To give an idea in terms of computation time, the hybrid algorithm took about 10 mn for the case with three random effects on a standard work station (Bi Xeon, 3.8 GHz).

\begin{table}
\caption{\label{meaniter}Percentage of convergence in less than 150 iterations and mean number of iterations to converge with different random effects for the Global Marquardt and Hybrid algorithm.}
\begin {center}
\begin{tabular}{c|cc|cc}
\hline
\multirow{2}{*}{Random effects} & \multicolumn{2}{c}{Global algorithm} &\multicolumn{2}{c}{Hybrid algorithm} \\
 &  nb iter  & \% success & nb iter  & \% success\\
\hline
$\tilde{\lambda}$                             & $35$   & $81\%$   & $11$   & $100\%$ \\
$\tilde{\lambda},~\tilde{\mu}_{T^*}$              & $54$   & $59\%$   & $17$   & $100\%$ \\
$\tilde{\lambda},~\tilde{\mu}_{T^*},~\tilde{\pi}$ & $71$  & $49\%$   & $25$   & $94\%$ \\ %
\hline
\end{tabular}

\end{center}
\end{table}%}
\subsection{Efficiency of the bias correction}

We estimated the bias of the corrected $\check \theta^{\btau}$ and uncorrected MHLE $\tetau$  using 500 replicas of a distribution with three random effects bearing on $\tilde{\lambda},~\tilde{\mu}_{T^*},~\tilde{\pi}$. We first examine the case where  $\btau =\btau^*=(0.2,0.2,0.2)$.  The biases of the uncorrected MHLE are of order $10^{-2}$ for all parameters. The correction reduces the biases to the order of $10^{-3}$ (except for one parameter), which seems negligible. 

\begin{table}[H]
\caption{\label{biais} Parameter values for the simulation and uncorrected and corrected biases}
\begin{center}
\begin{tabular}{c|c|cc|cc}
\hline
Parameters & True value & \multicolumn{2}{c}{Mean estimated value} &\multicolumn{2}{c}{Bias} \\
$$ & $$ & non corr  &corr & non corr  &corr\\
\hline
$\tilde \lambda$ & $4.10$   & $4.14$   & $4.10$   & $4.03~10^{-2}$ & $1.67~10^{-3}$ \\
$\tilde \mu_{T^*}  $ & $-1.60$ & $-1.54$ & $-1.60$ & $5.57~10^{-2}$ & $3.67~10^{-3}$ \\
$\tilde \pi    $ & $-0.170$   & $-0.160$   & $-0.166$   & $1.01~10^{-2}$ & $3.69~10^{-3}$ \\
$\gamma_0 $ & $-3.00$  & $-2.98 $ & $-3.00$  & $1.50~10^{-2}$ & $-3.60~10^{-3}$  \\
$\beta_1 $ & $-1.10$ & $-1.08$ & $-1.10$ & $2.04~10^{-2}$ & $2.59~10^{-3}$ \\
$\beta_2$  & $-1.40$ & $-1.35$ & $-1.39$ & $4.55~10^{-2}$ & $1.11~10^{-2}$ \\ %
\hline
\end{tabular}
\end{center}
\end{table}

\subsection{Property of the cMHLE}
We wished to check whether the asymptotic results hold in practice. We simulated data from the standard model of section \ref{standard}. In the first simulation (case 1) we took as standard deviations of the random effects $\tau^*_{\lambda}=\tau^*_{\mu_T}=\tau^*_{\mu_{T^*}}=0.2$. In a second simulation (case 2) we took $\tau^*_{\lambda}=0.1$, $\tau^*_{\mu_T}=0.2$, $\tau^*_{\mu_{T^*}}=0.3$. We did 500 replications and computed the root mean square errors (RMSE) and coverage rate of .95 confidence intervals of the estimated fixed parameters obtained in fixing the components of $\btau^u$ in the h-likelihood at values $\tau^u=0.1; 0.2; 0.3$. The results for the RMSE are shown in Table \ref{RMSE}. In the first case, the results tend to be better when $\tau^u=0.2$ which is closer to the $\tau^*_r$, while $\tau^u=0.3$ tends to be better than $\tau^u=0.1$. For the second case the results for $\tau^u=0.2$ and $\tau^u=0.3$ were approximately of the same quality, better than for $\tau^u=0.1$. It is striking that most of the RMSE are roughly of the same order, between $10^{-2}$ and $10^{-1}$. These RMSE can be interpreted as typical relative errors on the natural parameter; in these simulation the order of magnitude is about $5\%$. In term of coverage rates, the results (see Table \ref{coverage}) are not very good for $\tau^u=0.1$. They are satisfactory for  $\tau^u=0.2$ , and even more satisfactory for $\tau^u=0.3$. This corroborates our strategy based on a reasonable upper bound $\tau^u$ of the $\tau^*_r$.

%We see that the MHLE seem to be consistent and the variances correctly estimated whatever the value of $\btau$, although the variance is smaller, giving a better efficiency, for the value $0.2$  corresponding to that of the simulation.

\begin{table}[H]
\caption{Root Mean Square Error, $100$ subjects, $500$ replications. \label{RMSE}}
\begin {center}
\begin{tabular}{c|cc|cc|cc}
\hline
$\tau^u$ & \multicolumn{2}{c}{0.1} &\multicolumn{2}{c}{0.2}  & \multicolumn{2}{c}{0.3}\\
\hline
Par. & case 1  & case 2 & case 1  & case 2 & case 1  & case 2 \\
\hline
$\tilde \lambda$ & $5.62~10^{-2}$ & $4.47~10^{-2}$   & $3.34~10^{-2}$ & $3.94~10^{-2}$  & $4.95~10^{-2}$ & $3.76~10^{-2}$ \\
$\tilde \mu_{T^*}$&$9.10~10^{-2}$ & $7.67~10^{-2}$   & $4.64~10^{-2}$ & $6.84~10^{-2}$  & $8.50~10^{-2}$ & $7.47~10^{-2}$ \\
$\tilde \pi    $ & $7.27~10^{-2}$ & $5.95~10^{-2}$   & $5.14~10^{-2}$ & $5.86~10^{-2}$  & $5.25~10^{-2}$ & $5.71~10^{-2}$ \\
$\gamma_0 $ & $2.60~10^{-1}$ & $1.88~10^{-1}$   & $1.35~10^{-2}$ & $1.56~10^{-1}$  & $1.77~10^{-1}$ & $1.52~10^{-1}$ \\
$\beta_1$ & $1.74~10^{-1}$ & $1.59~10^{-1}$   & $1.01~10^{-1}$ & $1.04~10^{-1}$  & $1.04~10^{-1}$ & $9.89~10^{-2}$ \\
$\beta_2$ & $1.67~10^{-1}$ & $1.91~10^{-1}$   & $1.02~10^{-1}$ & $1.35~10^{-1}$  & $1.06~10^{-1}$ & $9.90~10^{-2}$ \\ %
\hline
\end{tabular}
\end{center}
\end{table}

%Finally in a third simulation we studied the choice of $\btau$ by crossvalidation. We simulated 500 replications with $\btau^*_{\lambda}=\btau^*_{\mu_T}=\btau^*_{\mu_{T^*}}=\btau=0.2$ as in the first simulation. We chose $\btau$ in the h-likelihood by the numerical CV and the approximated CV. For both methods we computed the mean and the standard deviation of the chosen $\btau$. Moreover we computed the correlation between a true future value of the viral load and the predicted value using either $\btau^*$ or the CV choices. The results are shown in Table 3...

\begin{table}[H]
\caption{Coverage rate, $100$ subjects, $500$ replications. \label{coverage}}
\begin {center}
\begin{tabular}{c|cc|cc|cc}
\hline
$\tau^u$ & \multicolumn{2}{c}{0.1} &\multicolumn{2}{c}{0.2}  & \multicolumn{2}{c}{0.3}\\
\hline
Par. & case 1  & case 2 & case 1  & case 2 & case 1  & case 2 \\
\hline
$\tilde \lambda$ & $90\%$ & $81\%$ & $95\%$ & $92\%$ & $89\%$ & $89\%$ \\
$\tilde \mu_{T^*}$&$94\%$ & $81\%$ & $94\%$ & $88\%$ & $93\%$ & $89\%$  \\
$\tilde \pi    $ & $96\%$ & $95\%$ & $94\%$ & $92\%$ & $93\%$ & $94\%$  \\
$\gamma_0 $ & $97\%$ & $97\%$ & $98\%$ & $96\%$ & $94\%$ & $95\%$  \\
$\beta_1$ & $94\%$ & $87\%$ & $97\%$ & $93\%$ & $93\%$ & $94\%$ \\
$\beta_2$ & $85\%$ & $74\%$ & $97\%$ & $93\%$ & $92\%$ & $93\%$ \\ %
\hline
\end{tabular}
\end{center}
\end{table}

\section{APPLICATION TO A CLINICAL TRIAL}
As an application of the proposed method, we aimed at estimating
the difference of treatment effects in a randomized clinical trial (Molina et al., 1999}). The ALBI ANRS 070
trial compared over 24 weeks the combination of zidovudine plus
lamivudine (AZT+3TC) with that of stavudine plus didanosine
(ddI+d4T) (a third arm alternating from one regimen to another was
not considered in this paper). The inclusion criteria were CD4 $\geq$ 200 cells$/\mu L$ and HIV RNA level between 4
and 5 $log_{10}$ copies/mL within 15 days before entry into the
study. The primary outcome measure {defined in the study protocol}
was the antiretroviral effect as measured by the mean change in
HIV RNA level between baseline and 24 weeks by use of the
ultra-sensitive PCR assay with lower limit of quantification of 50
copies/mL (1.7 $log_{10}$). In the main analysis of
Molina et al. (1999), HIV RNA values reported as $<$ 50 copies/mL were
considered equivalent to 50 copies/mL; 51 patients were included
in each treatment group. Over the 24-week period, HIV RNA level
declined in the two groups, with mean (SE) decreases at the end of
the study of 1.26 (0.09) $log_{10}$ copies/mL in the AZT+3TC group
and 2.26 (0.11) $log_{10}$ copies/mL in the ddI+d4T group.
% The mean increase in CD4 count was larger in ddI+d4T group than in AZT+3TC group (124 cells$/\mu L$ vs. 62 cells$/\mu L$, p=0.012).

 We used the model described in section 2.3. In this application only the first two components $Y_{ij1}$ and $Y_{ij2}$ were observed. Moreover only a left-censored version of $Y_{ij1}$ was observed; this was taken into account in the likelihood as in Guedj, Thi\'ebaut and Commenges (2007). In view of less informative observations than in the simulations we fixed the values of three parameters:  $\tilde \mu_{T}= -2.20$, $\tilde \mu_V =3.40$ and $\gamma_0= -3$. We put random effects on $\lambda$, $\pi$ and $\mu_T^*$, working with $\tau^u=0.3$. 
The estimated values of the parameters in natural logarithmic scale are displayed in Table 5. Reverting to natural parameters we find: $\hat \lambda= 56.8~ [54.1; 59.7]$; $\hat \pi=0.79~ [0.67; 0.92]$; $\hat \mu_{T^*}= 0.18~ [0.17; 0.19]$.
In addition it was possible to test whether the two treatment groups differed. The relevant null hypothesis is ``$\eta = 0$'', where $\eta=\beta_2-\beta_1$. A natural test statistic is $W={\hat \eta \over \sqrt{\widehat \var~ \hat \eta}}$, where $\hat \eta=\hat \beta_2-\hat \beta_1 $ and $\widehat \var ~\hat \eta$ can easily be computed from the estimate of the asymptotic variance matrix $\Sigma$. 

We found $\hat \eta= 0.242$, $\widehat \var ~\hat \eta = 5.16~10^{-3}$; this gives $W=3.37$ and a p-value equal to $p=7~10^{-4}$. Thus we conclude as expected that the treatment groups differ, and more precisely that the infectivity of the virus has been reduced more drastically in the ddI+d4T than in AZT+3TC group. Baseline infectivity is multiplied by a factor estimated to $e^{\hat \beta_2}=0.25$ and  $e^{\hat \beta_1}=0.32$ in the ddI+d4T than in AZT+3TC groups respectively.

\begin{table}[H]
\caption{Estimated parameters based on the ALBI clinical trial}
\begin{center}
\begin{tabular}{c|c|c|c}
\hline
Parameters & Uncorrected Values & Corrected Values & Confidence interval\\
\hline 
$\tilde \lambda$ & $4.05$ &$4.04$& $[3.99;4.09]$ \\
$\tilde \pi    $ & $-0.129 $ & $-0.242$ & $[-0.401;-0.083]$ \\
$\tilde \mu_{T^*}    $ & $-1.74$ & $-1.73$ & $[-1.80;-1.65]$ \\
$\beta_1$ & $-1.33$ & $-1.12$& $[-1.29;-0.957]$ \\
$\beta_2$ & $-1.53$ & $-1.37$& $[-1.56;-1.17]$ \\
$\sigma_{CD4}$ & $0.173$ &$0.168$ & $[0.151;0.185]$ \\
$\sigma_{CV}$ & $0.584$ & $0.541$ & $[0.501;0.582]$\\  %
\hline
\end{tabular}
\end{center}
\end{table}

\section{CONCLUSION}
We have developed a hierarchical likelihood approach for inference in an HIV dynamical model. We have obtained the asymptotic distribution of the MHLE, we have derived a procedure which makes the bias negligible and we have developed an efficient maximization  algorithm. Our simulations show that the whole approach works.

We have shown that it could be applied to the analysis of a real data set. Rather precise estimates of the parameters were obtained. One limitation of this approach is that some parameters must be fixed because of identifiability problems. The model itself, although it is already statistically challenging, may be too simple from a biological point of view. The development of such an approach would require richer data, for instance observing the number of infected T cells. 

The main advantage of this approach is that it is easy to implement and very fast as compared to the two main competing approaches, likelihood and Bayesian inference. The main limitation is that it does not attempt to estimate the variances of the random effects. In our application we already have a knowledge of the range of values of these variances. Thus the method can be used for exploring possible models while likelihood or Bayesian inference can be used when estimates of the variances of the random effects are needed.\vspace{2mm}

\noindent {\bf Acknowledgments.} The authors thank the investigators of the ALBI ANRS-070 clinical trial and particularly J. M. Molina (principal investigator) and G. Ch\^ene (methodologist).

\section*{REFERENCES}
\noindent
\setlength{\parindent}{-8mm}

Beal, S.L. and Sheiner, L.B. (1982) Estimating population kinetics. {\em Critical Reviews in Biomedical Engineering}, {\bf 8}, 195-222.%\vspace{2mm}

Davidian, M. and Giltinan, D.M. (1995) {\em Nonlinear models for repeated measurements data}, Chapman \& Hall.%\vspace{2mm}

Donnet, S. and Samson, A. (2007) Estimation of parameters in incomplete data models defined by dynamical systems. {\em Journal of Statistical Planning and Inference}, {\bf 137}, 2815-2831%\vspace{2mm}

Bradley Efron, R.J. Tibshirani (1993) {\em An Introduction to the Bootstrap.} Chapman \& Hall.%\vspace{2mm}

Fletcher, R. (1987) {\em Practical Methods of Optimization}. John Wiley \& Sons (Chichester).%\vspace{2mm}

Guedj, J., Thi\'ebaut, R. and Commenges, D. (2007) Maximum likelihood estimation in dynamical models of HIV. {\em Biometrics}, {\bf 63}, 1198-1206.%\vspace{2mm}

Ho, D.D., Neumann, A.U., Perelson, A.S., Chen, W., Leonard, J.M. and Markowitz, M.(1995) Rapid turnover of plasma virions and CD4 lymphocytes in HIV-1 infection. {\em Nature}, {\bf 373}, 123-126.%\vspace{2mm}

Huang, X., Liu, D and Wu, H. (2006) Hierarchical Bayesian methods for estimation of parameters in a longitudinal HIV dynamic system. {\em Biometrics}, {\bf 62}, 413,423.%\vspace{2mm}

Kuhn, E. and Lavielle, M. (2005) Maximum likelihood estimation in nonlinear mixed effects models. {\em Computational Statistics \& Data Analysis}, {\bf 49}, 1020-1038.%\vspace{2mm}

Lee, Y. and Nelder, J.A. (1996) Hierarchical Generalized Linear Models.
{\em Journal of the Royal Statistical Society. Series B}, {\bf 58}, pp. 619-678%\vspace{2mm}

Lee, Y. and Nelder, J.A. (2001) Hierarchical generalised linear models: A synthesis of generalised linear models, random-effect models and structured dispersions. {\em Biometrika}, {\bf 88}, 987-1006%\vspace{2mm}

Lee, Y., Nelder, J.A. and Pawitan, Y. (2006) {\em Generalized linear models with random effects}, Chapman and Hall.%\vspace{2mm}

Lindstrom, M. and  Bates, D. (1990) Nonlinear mixed effects models for repeated measures data. {\em Biometrics}, {\bf 46}, 673-687.%\vspace{2mm}

McGilchrist, C.A. and Aisbett,  C.W. (1991) Regression with Frailty in Survival Analysis.  {\em Biometrics}, {\bf 47}, 461-466.%\vspace{2mm}

Marquardt, D. (1963) An algorithm for least-squares estimation of nonlinear parameters. {\em SIAM Journal of Applied Mathematics}, {\bf 11}, 431-441.%\vspace{2mm}

Molina, J.M., Ch\^ene, G., Ferchal, F., Journot, V., Pellegrin, I., Sombardier, M. N.,  Rancinan, C., Cotte, L., Madelaine, I., Debord, T. and Decazes, J. M. (1999) The {ALBI} trial: A Randomized Controlled Trial Comparing Stavudine Plus Didanosine with Zidovudine Plus Lamivudine and a Regimen Alternating Both Combinations in Previously untreated Patients Infected with Human immunodeficiency Virus. {\em The Journal of Infectious Diseases}, {\bf 180}, 351-358.%\vspace{2mm}

Noh, M. and Lee, Y. (2008) Hierarchical-likelihood approach for nonlinear mixed-effects models. {\em Computational Statistics \& Data Analysis}, {\bf 52}, 3517-3527.%\vspace{2mm}

Nowak, M.A. and Bangham, C.R.M. (1996) Population dynamics and immune response to persistent viruses. {\em Science}, {\bf 272}, 74-79.%\vspace{2mm}

Nowak, M.A. and  May R.M. (2000) {\em Virus Dynamics: Mathematical Principles of Immunology and Virology}, Oxford University Press.%\vspace{2mm}

O'Sullivan, F. (1988) Fast computation of fully automated log-density and log-hazard estimators. {\em SIAM Journal on Scientific and Statistical Computing}, {\bf 9}, 363-379.%\vspace{2mm}

Putter, H, Heisterkamp, S.H., Lange, J.M. and de Wolf, F. (2002) A Bayesian approach to parameter estimation in HIV dynamical models. {\em Statistics in Medicine}, {\bf 21}, 2199-2214.%\vspace{2mm}

Ramsay, J. O., Hooker, G., Campbell, D. and Cao, J. (2007) Parameter estimation for differential equtions: a generalized smoothing approach. {\em Journal of the Royal Statistical Society: Series B}, {\bf 69}, 741-796.

Therneau, T.M. and Grambsch P.M. (2000) {\em Modeling survival data: extending the Cox model}, Springer.%\vspace{2mm}

van der Vaart, A. (1998) {\em Asymptotic Statistics}, Cambridge.%\vspace{2mm}

Wu, H. (2005) Statistical methods for HIV dynamic studies in
AIDS clinical trials. {\em Statistical Methods in Medical Research}, {\bf 14}, 171-192.

Wu, H. and Ding, A. (1999) Population HIV-1 Dynamics in Vivo: Applicable Models and Inferential Tools for Virological Data from AIDS Clinical Trials. {\em Biometrics}, {\bf  55}, 410-418.
\end{document}